\documentclass{ifacconf}

\usepackage{graphicx}      %
\usepackage{natbib}        %
\usepackage[table]{xcolor}
\usepackage{amsmath}
\usepackage{tabularx,booktabs}
\usepackage[export]{adjustbox}
\allowdisplaybreaks
\usepackage{hhline}
\usepackage{url}
\pdfoutput=1

\definecolor{MyBlue}{HTML}{A9CCE3}
\newtheorem{remark}{\textbf{Remark}}

\makeatletter
\let\old@ssect\@ssect %
\makeatother

\usepackage{natbib}
\usepackage{hyperref}
\hypersetup{
	colorlinks=true,
	linkcolor=blue,
	filecolor=black,      
	urlcolor=cyan,
	citecolor=black
}

\makeatletter
\def\@ssect#1#2#3#4#5#6{%
	\NR@gettitle{#6}%
	\old@ssect{#1}{#2}{#3}{#4}{#5}{#6}%
}
\makeatother

\begin{document}
\begin{frontmatter}

\title{Structured Online Learning-based Control of Continuous-time Nonlinear Systems  \thanksref{footnoteinfo}} 

\thanks[footnoteinfo]{This work was supported in part by an NSERC Discovery Grant, the Canada Research Chairs program, and an Early Researcher Award from the Ontario Ministry of Research, Innovation and Science.}

\author[First]{Milad Farsi} 
\author[First]{Jun Liu}

\address[First]{Applied Mathematics Department, University of Waterloo, Canada (e-mail: mfarsi@uwaterloo.ca, j.liu@uwaterloo.ca)}

\begin{abstract}                %
 Model-based reinforcement learning techniques accelerate the learning task by employing a transition model to make predictions. In this paper, a model-based learning approach is presented that iteratively computes the optimal value function based on the most recent update of the model. Assuming a structured continuous-time model of the system in terms of a set of bases, we formulate an infinite horizon optimal control problem addressing a given control objective. The structure of the system along with a value function parameterized in the quadratic form provides a flexibility in analytically calculating an update rule for the parameters. Hence, a matrix differential equation of the parameters is obtained, where the solution is used to characterize the optimal feedback control in terms of the bases, at any time step. Moreover, the quadratic form of the value function suggests a compact way of updating the parameters that considerably decreases the computational complexity. Considering the state-dependency of the differential equation, we exploit the obtained framework as an online learning-based algorithm. In the numerical results, the presented algorithm is implemented on four nonlinear benchmark examples, where the regulation problem is successfully solved while an identified model of the system is obtained with a bounded prediction error.   

\end{abstract}


\end{frontmatter}

\section{Introduction}
Model-based reinforcement learning techniques, as apposed to direct methods in learning, are known to be more data-efficient. Direct methods usually require enormous data and hours of training even for simple applications \citep{SOL:duan2016benchmarking}, while model-based techniques can show optimal behavior in a limited number of trials. This property, in addition to the flexibilities in changing learning objectives and performing farther safety analysis make them more suitable for real-world implementations, such as robotics. In model-based approaches, having a deterministic or probabilistic description of the transition system saves much of the effort spent by direct methods in treating any point in the state-control space individually.  Hence, the role of model-based techniques become even more significant when it comes to problems with continuous control rather than discrete actions (\cite{SOL:sutton1990integrated,SOL:atkeson1997comparison,SOL:powell2004handbook}). 

Various model-based approaches can be found in the literature that are mainly categorized under two topics: value function and policy search methods. In value methods, known also as approximate/adaptive dynamic programming techniques, a value function is used to construct the policy. However, policy search methods directly improve the policy to achieve optimality. A review of recent techniques can be found in \cite{SOL:2017reinforcement,SOL:recht2019tour,SOL:polydoros2017survey,SOL:kamalapurkar2018model}. Another survey is given in \cite{SOL:benosman2018model} that contrasts model-based techniques with data-driven adaptive techniques.

Value methods in reinforcement learning normally require solving the well-known Hamilton-Jacobi-Bellman (HJB). However, common techniques for solving such equations suffer from curse of dimensionality. Hence, in approximate dynamic programming techniques, a parametric or non-parametric model is used to approximate the solution. In \cite{SOL:lewis2009reinforcement}, some related approaches are reviewed that fundamentally follow the actor-critic structure (\cite{SOL:barto1983neuronlike}), such as value and policy iteration algorithms.

In such approaches, the Bellman error, which is obtained from the exploration of the state space, is used to improve the parameters estimated in a gradient-descent or least-square loop that require persistent excitation condition. Since the Bellman error obtained is only valid along the trajectories of the system, sufficient exploration in the state-space is required to efficiently estimate the parameters. \cite{SOL:kamalapurkar2018model} has reviewed different strategies employed to increase the data-efficiency in exploration. In \cite{SOL:vamvoudakis2012multi,SOL:modares2014integral}, a probing signal is added to the control to enhance the exploring properties of the policy. In another approach, the recorded data of explorations is used as a replay of experience to increase the data efficiency. In \cite{,SOL:modares2014integral}, the model obtained from identification is used to acquire more experience by doing simulation in an offline routine that decreases the need for visiting any point in the state space.

As an alternative method, considering a nonlinear control affine system with a known input coupling function, \cite{SOL:kamalapurkar2016model} used a parametric model to approximate the value function. Then, they employed a least-square minimization technique to adjust the parameters according to the Bellman error which can be calculated at any arbitrary point of the state space by having identified internal dynamics of the system and approximated state derivatives. 

In this paper, as a different approach for obtaining the parameters of value function, we aim on solving a closed-loop optimal control problem online, rather than using a minimization technique to update the parameters based on the Bellman error. Although, Bellman theory is clearly behind all these techniques, the formulation here is contrasted with previous approaches: assuming a particular structure for the identified system allows us to analytically obtain a matrix differential equation in terms of parameters. This relaxes the need of a gradient-descent or least-square technique that is used in similar approaches. Unlike \cite{SOL:zhang2011data,SOL:bhasin2013novel,SOL:kamalapurkar2016model}, we assume a quadratic form for the value function that provides a compact way for parameterizing the value in terms of quadratic terms, following the structure considered for the dynamics. Accordingly, we will refer the presented framework as a structured online learning (SOL) algorithm. 

In \cite{SOL:brunton2016discovering}, an algorithm for sparse identification of nonlinear dynamical systems (SINDy) is presented to obtain the explicit dynamics of the system. In \cite{SOL:kaiser2018sparse}, a control approach is employed based on SINDy that includes two independent stages: the identification with a generated random signal and the model predictive control. This technique is shown to be a data-efficient identification scheme that, in addition, can handle the noise in data. Although the identification methods can be used in SOL algorithm are not limited to any particular approach, as another contribution of this paper, we employ SINDy to achieve faster convergence of the parameters in an online scheme. This is contrasted with \cite{SOL:kaiser2018sparse} in a way that we iteratively perform both identification and control in a single loop.    

The rest of the paper is organized as follows. In Section \ref{two}, we propose an optimal control approach based on a particular structure of dynamics, and characterize the optimal feedback control based on a matrix of parameters obtained by a differential equation. Section \ref{three} outlines the SOL algorithm designed based on the obtained results. In Section \ref{four}, we present the numerical results of this algorithm implemented on a few benchmark examples.  
\subsection{Notations}
We will denote $p$-norm by $\|\cdot\|_p$. For defining a set of basis, any operator on vector $x$, is performed component-wise e.g. $x^2=[x_1^2,\dots,x_n^2]$. Moreover, a diagonal square matrix $A$ with elements $A_1,\dots,A_n$ on the diagonal is shortened as $A=\text{diag}([A_1,\dots,A_n])$.

\section{A structured Approximate Optimal Control Framework }{
	\label{two}
	Consider the nonlinear affine system
	\begin{align} \label{sys}
	\dot{x}=f(x)+g(x)u
	\end{align}
	where $x\in D\subset {\rm I\!R}^n$, $u\in\Omega\subset{\rm I\!R}^m$, and $f,g:D\rightarrow {\rm I\!R}^n$.
	
	The cost function to be minimized along the trajectory, started from the initial condition $x_0=x(0)$, is considered in the following linear quadratic form
	\begin{align} \label{cost_x}
	J(x_0,u)=\lim_{T\rightarrow\infty}\int_{0}^T \e^{-\gamma t}\left(x^TQx+u^TRu\right)  dt,
	\end{align}
	where $Q\in{\rm I\!R}^{n\times n}$ is positive semi-definite, $\gamma\geq 0$ is the discount factor, and $R\in{\rm I\!R}^{m\times m}$ is a diagonal matrix with only positive values, given by design criteria.
	
	For the closed-loop system, by assuming a feedback control law $u=\omega(x(t))$ for $t\in [0,\infty)$, the optimal control is given by 
	\begin{align}\label{2661}
	\omega^*= \quad \arg\min_{\textcolor{black}{u(\cdot) \in \Gamma(x_0)}} J(x_0, u(\cdot)),
	\end{align} 
	where $\Gamma$ is the set of admissible controls. 
	\subsubsection{Assumption 1.}\label{assumption} {
		$f$ and $g$ can be identified or effectively approximated within the domain of interest by the linear combination of some basis functions $\phi_i\in C^1:D\rightarrow {\rm I\!R}$ for $i=1,2,\dots,p$.
	
	Accordingly, (\ref{sys}) is rewritten as	
	\begin{align} \label{sys_approx}
	\dot{x}=W\Phi(x)+\sum_{j=1}^{m} W_j\Phi(x)u_j,
	\end{align}
	where $W$ and $W_j\in{\rm I\!R}^{n\times p}$ are the matrices of the coefficients obtained for $j=1,2,\dots,m$, and $\Phi(x)=[ \phi_{1}(x) \quad\dots\quad \phi_{p}(x)]^T$. 
	
	In what follows, without loss of generality, the cost defined in (\ref{cost_x}) is transformed to the space of bases $\Phi(x)$, that is
		\begin{align}
	J(x_0,u)=\lim_{T\rightarrow\infty}\int_{0}^T \e^{-\gamma t}\left(\Phi(x)^T\bar{Q}\Phi(x)+u^TRu \right) dt,
	\end{align}
	where $\bar{Q}=\text{diag}\left([Q],[\boldsymbol 0_{(p-n)\times (p-n)}]\right)$ is a block diagonal matrix that contains all zeros except the first block Q which correspond to the linear bases $x$.
	
	Then the corresponding HJB equation can be written by the  Hamiltonian defined as
	\begin{align}\label{HJB}
	-&{\frac{\partial}{\partial{t}}}({\e^{-\gamma t}V})\hspace{-1mm}=\hspace{-1mm}\min_{{u(\cdot) \in \Gamma(x_0)}}\{H=\e^{-\gamma t}\big(\Phi(x)^T\bar{Q}\Phi(x)+u^TRu\big)\nonumber\\ 
	&\qquad \qquad\qquad{+\e^{-\gamma t}\frac{\partial{V}}{\partial{x}}}^T(W\Phi(x)+\sum_{j=1}^{m} W_j\Phi(x)u_j)\} .
	\end{align}
	In general, there exist no analytical approach that can solve such partial differential equation and obtain the optimal value function. However, it has been shown in the literature that approximate solutions can be computed by numerical techniques. 
	
	Assume the optimal value function in the following form
	\begin{align}\label{value}
	V=\Phi(x)^TP\Phi(x),
	\end{align}
	where $P$ is symmetric.
	\begin{remark}
		Unlike other approximate optimal approaches in the literature, such as, \cite{SOL:zhang2011data,SOL:bhasin2013novel,SOL:kamalapurkar2016efficient}, that use a linear combination of bases to parameterize the value function, we assume a quadratic form. As a result, the value function now is defined in the product space $\Lambda:=\Phi \times\Phi$. Hence, it is expected that the resulting quadratic terms better contribute to basing a positive value function around $x=0$. Furthermore, due to the function-approximating properties of bases $\Phi$ itself, one may bring them to $\Lambda$ in addition by including a constant basis $c$ in $\Phi$. Therefore, compared to other approaches, the structure used in (\ref{value}) suggests a more compact way of formulating the problem where by only a limited number of bases in $\Phi$, we can attain a richer set $\Lambda$ to parameterize the value function.   
    \end{remark}
		Then the Hamiltonian is given by
	\begin{align}
	H=&\e^{-\gamma t}(\Phi(x)^T\bar{Q}\Phi(x)+u^TRu)\nonumber\\
	&+\e^{-\gamma t}\Phi(x)^TP\frac{\partial{\Phi(x)}}{\partial{x}}\bigg (W\Phi(x)+\sum_{j=1}^{m} W_j\Phi(x)u_j\bigg)\nonumber\\ \nonumber
	+&\e^{-\gamma t}\bigg(\Phi(x)^TW^T+\sum_{j=1}^{m} u_j^T\Phi(x)^TW_j^T\bigg )\frac{\partial{\Phi(x)}}{\partial{x}}^TP\Phi(x).\nonumber
	\end{align}

	Moreover, based on the structure of $R$, the quadratic term of $u$ is rewritten in terms of its components.
	\begin{align}
	&H=\nonumber\\
	&\e^{-\gamma t}\Bigg(\Phi(x)^T\bar{Q}\Phi(x)+\sum_{j=1}^{m}r_ju_j^2+\Phi(x)^TP\frac{\partial{\Phi(x)}}{\partial{x}}W\Phi(x)+\nonumber\\
	&\Phi(x)^TP\frac{\partial{\Phi(x)}}{\partial{x}}\bigg(\sum_{j=1}^{m} W_j\Phi(x)u_j\bigg)\hspace{-1.2mm}+\hspace{-0.8mm}\Phi(x)^TW^T\frac{\partial{\Phi(x)}}{\partial{x}}^TP\Phi(x) \nonumber\\ 
	&+\bigg(\sum_{j=1}^{m} u_j\Phi(x)^TW_j^T\bigg)\frac{\partial{\Phi(x)}}{\partial{x}}^TP\Phi(x)\Bigg),
	\end{align}	
	where $r_j\neq 0$ is the $j$th component on the diagonal of matrix $R$.
	To minimize the resulting Hamiltonian we need
	\begin{align}
	\frac{\partial{H}}{\partial{u_j}}&=2r_ju_j+2\Phi(x)^TP\frac{\partial{\Phi(x)}}{\partial{x}}W_j\Phi(x)\\ \nonumber
	&=0, \quad \quad \quad j=1,2,\dots,m.
	\end{align}
	Hence, the $j$th optimal control input is obtained as
	\begin{align} \label{control}
	u_j^*=-\Phi(x)^Tr_j^{-1}P\frac{\partial{\Phi(x)}}{\partial{x}}W_j\Phi(x).
	\end{align} 
	By plugging in the optimal control and the value function in (\ref{HJB}) we get
	\begin{align}
	-&\e^{-\gamma t}\Phi(x)^T\dot{P}\Phi(x)+\gamma\e^{-\gamma t}\Phi(x)^TP\Phi(x)\nonumber\\ 
	&=\e^{-\gamma t}\Bigg(\Phi(x)^T\bar{Q}\Phi(x)+\nonumber\\ \nonumber
	&+\Phi(x)^TP\frac{\partial{\Phi(x)}}{\partial{x}}\bigg(\sum_{j=1}^{m}W_j\Phi(x)r_j^{-1}\Phi(x)^TW_j^T\bigg)\\ \nonumber &\hspace{60mm}\frac{\partial{\Phi(x)}}{\partial{x}}^TP\Phi(x) \\ \nonumber
	&-2\Phi(x)^TP\frac{\partial{\Phi(x)}}{\partial{x}}\bigg(\sum_{j=1}^{m} W_j\Phi(x)r_j^{-1}\Phi(x)^TW_j^T\bigg)\\ \nonumber &\hspace{60mm}\frac{\partial{\Phi(x)}}{\partial{x}}^TP\Phi(x)\\ \nonumber
	&+\Phi(x)^TP\frac{\partial{\Phi(x)}}{\partial{x}}W\Phi(x)+\Phi(x)^TW^T\frac{\partial{\Phi(x)}}{\partial{x}}^TP\Phi(x)\Bigg). \nonumber
	\end{align}
	This is rewritten as
\begin{align}
	&\Phi(x)^T\dot{P}\Phi(x)+\gamma\Phi(x)^TP\Phi(x)=\Phi(x)^T\bar{Q}\Phi(x)+\nonumber\\ \nonumber
	&-\Phi(x)^TP\frac{\partial{\Phi(x)}}{\partial{x}}\bigg(\sum_{j=1}^{m}W_j\Phi(x)r_j^{-1}\Phi(x)^TW_j^T\bigg)\\ \nonumber &\hspace{60mm}\frac{\partial{\Phi(x)}}{\partial{x}}^TP\Phi(x) \\ \nonumber
	&+\Phi(x)^TP\frac{\partial{\Phi(x)}}{\partial{x}}W\Phi(x)+\Phi(x)^TW^T\frac{\partial{\Phi(x)}}{\partial{x}}^TP\Phi(x),
	\end{align} 
	where a sufficient condition to hold this equation is 
	\begin{align} \label{P_ode}
	-\dot{P}=&\bar{Q}+P\frac{\partial{\Phi(x)}}{\partial{x}}W+W^T\frac{\partial{\Phi(x)}}{\partial{x}}^TP-\gamma P\nonumber\\ 
	&-P\frac{\partial{\Phi(x)}}{\partial{x}}\bigg(\sum_{j=1}^{m}W_j\Phi(x)r_j^{-1}\Phi(x)^TW_j^T\bigg)\frac{\partial{\Phi(x)}}{\partial{x}}^TP.
	\end{align}
	This equation has to be solved backward to get a value of $P$ that characterizes the optimal value function (\ref{value}) and control (\ref{control}), however it has been shown that the forward integration of such equation converge to similar results as long as we are not very close to the initial time.

	\begin{remark} 
		While the similarity between the derived optimal control and the linear quadratic regulation (LQR) problem cannot be denied, there are substantial differences. It should be noted that the matrix differential equation (\ref{P_ode}) is derived in terms of $\Phi$ which is of dimension $p$ in contrast with the LQR formulation that includes only the linear terms of the state with dimension $n$. As a result, by using the proposed SOL framework, once we reached to the neighborhood of the goal ($x\equiv 0$), we have a richer value function that can be generalized to a larger domain than LQR which is only valid in the neighborhood of the origin. Fig.~\ref{fig:pend} and Fig. ~\ref{fig:pend1} denote the comparison results with the LQR control on a pendulum, where we choose $Q$ in a way that both states are penalized equally. Although they illustrate similar responses around the equilibrium points, LQR fails to effectively penalize the angular velocity when $x_0$ is not in the neighborhood of the equilibrium point.
	\end{remark}
	\begin{figure}[h]
	\begin{center}
		\hspace*{-5mm}  
		\includegraphics[width=7.4cm]{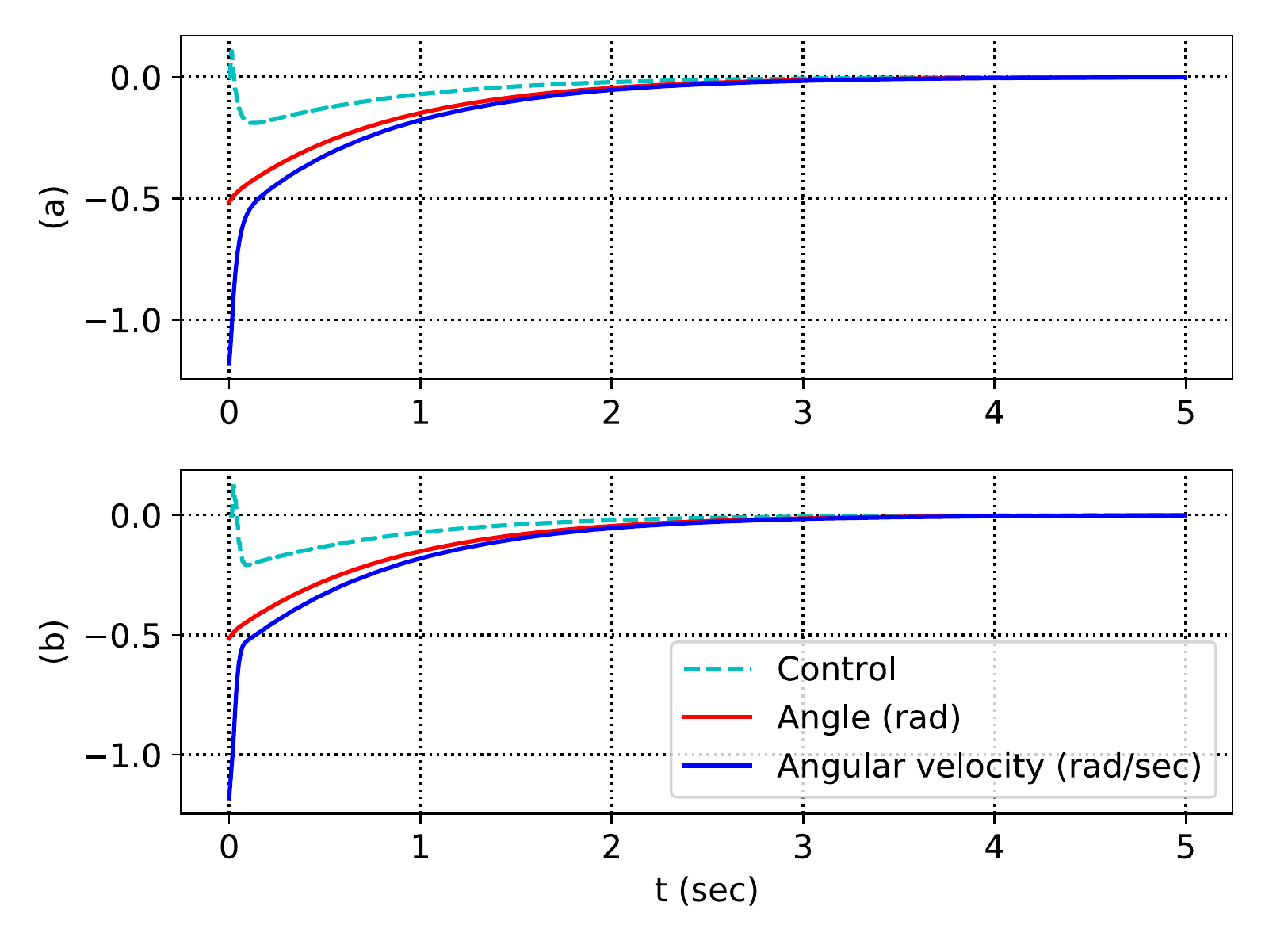}    %
		\caption{Comparison results of the pendulum for an initial condition ($x_0=[-0.51,-1.18]$) around the equilibrium point that shows very similar results . (a) LQR control (b) Control obtained by SOL. } 
		\label{fig:pend}
	\end{center}
\end{figure}

\begin{figure}[h]
	\begin{center}
		\includegraphics[width=7cm]{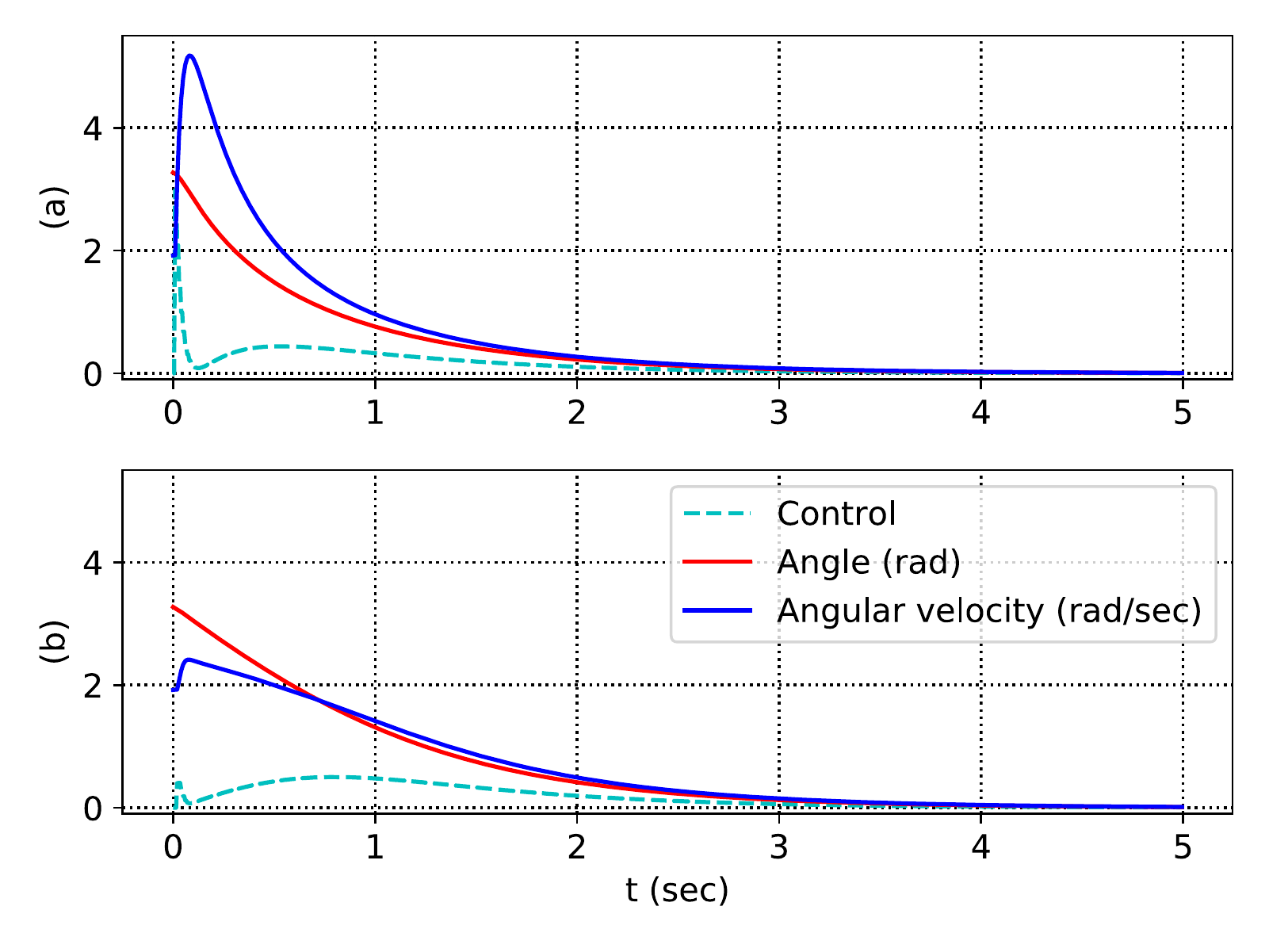}    %
		\caption{Comparison results of the pendulum for an initial condition ($x_0=[3.27,1.92]$) not in the neighborhood of the equilibrium point that shows a better penalization of the responses for SOL. (a) LQR control (b) Control obtained by SOL. } 
		\label{fig:pend1}
	\end{center}
\end{figure}
	\begin{remark}
		 In fact, it is not difficult to prove the LQR control as a special case of the proposed framework by limiting the bases to $\Phi(x)=[ (x_1\quad \dots \quad x_n) \quad 1]$. In this case, $P$ will be a $(2\times 2)$-block-structured matrix where the $\dot P_{11}$ block of (\ref{P_ode}) produces the well-known Riccati equation, with $\gamma=0$.
	\end{remark}
	\begin{remark}
		Because of the general case considered in obtaining (\ref{P_ode}), where $\Phi$ includes arbitrary basis functions of the state, there exists no way to escape from the state-dependency in this equation, except in the linear case as mentioned. Hence, we require (\ref{P_ode}) be solved along the trajectories of the system.    	
	\end{remark}
	
	In the next section, we will establish an online learning algorithm based on the proposed optimal control framework. 
}

\section{Structured online learning (SOL) algorithm}
\label{three}
By considering a general description of the nonlinear input affine  system in terms of some bases as in (\ref{sys_approx}), we obtained a structured optimal control framework that suggests using the state-dependent matrix differential equation (\ref{P_ode}) to achieve the parameters of the nonlinear feedback control. Next, we exploit this framework to propose the SOL algorithm. Hence the focus of this section will be on the algorithm and practical properties of SOL. 

The learning procedure is done in the following order: We first initialize $P$ with a zero matrix. In the control loop, at any time step $t_k$, the samples of the states are acquired and the set of bases are evaluated accordingly. Next, by an identification technique, the system model is updated. The measurements and the most recent model coefficients are used to integrate (\ref{P_ode}) and improve $P$, which is required to calculate the control value for the next step $t_{k+1}$. In what follows, we discuss the steps involved in more details with focusing on SINDy algorithm.       

\subsubsection{ODE Solver and Control Update:}In this approach, we run the system from some $x_0\in D$, then solve the matrix differential equation (\ref{P_ode}) along the trajectories of the system. Different solvers are already developed that can efficiently integrate differential equations. In the simulation, we use a Runge–Kutta solver to integrate the dynamics of the system that replaces the real system in a real-world application. Although the solver may take smaller steps, we only allow the measurements and control update at time steps $t_k=kh$, where $h$ is the sampling time and $k=0,1,2,\dots$ . For solving (\ref{P_ode}) in continuous time, we use the Runge–Kutta solver with a similar setting, where the weights and the states in this equation are updated by a system identification algorithm and the measurements $x_k$ at each ieteration of the control loop, respectively. A recommended choice for $P_0$ is a matrix with components of zero or very small values.

The differential equation (\ref{P_ode}) also requires evaluations of $\partial \Phi/\partial x_k$ at any time step. Since the bases $\Phi$ are chosen beforehand, the partial derivatives can be analytically calculated and stored as functions. Hence, they can be evaluated for any $x_k$ in a similar way as $\Phi$ itself. By solving (\ref{P_ode}), we can calculate the control update at any time step $t_k$ according to  (\ref{control}). Although, at the very first steps of learning, control is not expected to take effective steps toward the control objective, it can help in exploration of the state space and gradually improve by learning more about the dynamics. 

\begin{remark}
	The computational complexity of updating parameters by relation (\ref{P_ode}) is bounded by the complexity of matrix multiplications of dimension $p$ which is $\mathcal{O}(p^3)$. Moreover, it should be noted that, regarding the symmetry in the matrix of parameters $P$, this equation updates $L=(p^2+p)/2$ number of parameters which correspond to the number of bases used in the value function. Therefore, in terms of the number of parameters, the complexity of the proposed technique is $\mathcal{O}(L^{3/2})$. However, for instance, if 
	recursive least squares technique were employed with the same number of parameters, the computations are bounded by ${O}(L^{3})$. As a result, the proposed parameter update scheme can be done considerably faster than similar model-based techniques, such as, \cite{SOL:kamalapurkar2016model,SOL:bhasin2013novel}. In another effort, \cite{SOL:kamalapurkar2016efficient} decreased the number of bases used to improve the computational efficiency, while the complexity still remained as ${O}(L^{3})$.             
\end{remark} 
\subsubsection{Identified Model Update:}We considered a given structured nonlinear system as in Assumption 1. Therefore, having the control and state samples of the system, we need an algorithm that updates the weights in (\ref{sys_approx}) corresponding to a given cost $E(\cdot)$ as
\begin{align}\label{learning_cost}
{[W\quad W_1\dots W_j]}_k= \arg\min_{\hat W_k} \quad E(\dot x_k, \hat W_k, \Theta(x_k,u_k) ),
\end{align} 
where $\Theta(x_k,u_k)={[\Phi^T(x_k)\quad \Phi^T(x)u_1\quad \dots \quad \Phi^T(x)u_m]}^T_k$, and $k$ is the time step.

 As studied in \cite{SOL:brunton2016discovering,SOL:kaiser2018sparse}, SINDy is a data-efficient tool to extract the underlying sparse dynamics of the sampled data. Hence, we use SINDy to update the weights of the system to be learned. In this approach, along with the identification, the sparsity is also promoted in the weights by minimizing 
 \begin{align}
 E(\dot x, \hat W, \Theta(x_k,u_k) )=\|\dot x_k-\hat W_k\Theta(x_k,u_k)\|^2_2+\lambda\|\hat W_k\|_1 \nonumber
 \end{align}
 as in (\ref{learning_cost}), where $\lambda$ is a positive constant. Furthermore, there exist other techniques that can be alternatively used, such as neural networks, nonlinear regression, or any other function approximation and system identification methods that can minimize the prediction error $\|\dot x_k-\hat W_k\Theta(x_k,u_k)\|$ by assuming the structure defined in (\ref{sys_approx}). 

\subsubsection{Database Update:} For using SINDy algorithm, a database of samples is required to recursively perform regressions. These weights correspond to a library of functions given in $\Phi$. Any sample of the system, to be stored in the database at time $k$, includes $\Theta(x_k,u_k)$ and the derivatives of the states approximated by $(x_{k}-x_{k-1})/h$.

We adopt SINDy to do an online learning task, meaning that the database has to gradually build along with the exploration and control. Different approached can be employed in the choice of samples and building a database online. A comparison of these techniques can be found in in \cite{SOL:kivinen2004online,SOL:van2014online}. In the implementations of SOL done in this paper, we assume a given maximum size of database $N_d$, then we keep adding the samples with larger prediction errors to the database. This is done in a loop together with updating the control until the bound of the prediction error obtained allows the control to regulate the system to the given reference state. If the maximum number of samples in database is reached, we forget the oldest sample and replace it with the recent one.     

\section{Simulation Results}
\label{four}
We have implemented the proposed approach on four examples which, are presented in two categories, considering Assumption 1: 1) the dynamics can be written exactly in terms of some choice of basis functions, and 2) the dynamics includes some terms that are required to be approximated in the space of some given bases. 

As mentioned, in these numerical examples we have exploited the SINDy algorithm for the identification purpose, however, clearly the focus of the simulations here is on the properties of the proposed control scheme rather than the identification part, regarding that SINDy already has been extensively studied in \cite{SOL:brunton2016discovering,SOL:kaiser2018sparse} as an offline identification algorithm. The SINDy algorithm adopted here is a powerful tool to obtain the dynamics of the system with a good precision. However, this depends greatly on how efficiently we can approximate the derivatives of the states. Hence, in different implementations, higher sampling rates or higher order approximation of the derivatives may be needed. For the same reason, in the proposed examples, the number of samples used and the system obtained may be further tuned to match the level of quality reported in \cite{SOL:brunton2016discovering,SOL:kaiser2018sparse}.

The simulations are done in Python, where we used Vpython module (\cite{SOL:scherer2000vpython}) to generate the graphics. We have set the sampling rate to $200$Hz ($h=5$ms) for all the examples, unless explicitly mentioned otherwise. The control input value is updated at every other time step meaning that the update rate is $100$Hz. The simulation is stopped if the trajectory reaches to the boundary of $D$ or a timeout is reached without satisfying the objective. Moreover, if the regulation objective is to reach a point other than the origin ($x\equiv0$), we consider the cost (\ref{cost_x}), the value (\ref{value}), and the obtained differential equation of the value parameters (\ref{P_ode}) by redefining $x:=x-x_\text{ref}$.     
\subsection{Systems Identifiable in Terms of a Given Set of Bases}
In the following two examples, we assume that the bases constituting the systems dynamics exist in $\Phi$. The system identified, after running the proposed learning algorithm and obtaining the value function, clearly depends on the identification algorithm used and its tuning parameters.

 In Table~\ref{tb:pend_lorenz}, we illustrate the variations of the identified system and the corresponding value function by implementing the presented SOL algorithm with the exact $\dot x$ and with the first order approximation of the derivative. It can be observed that, in the pendulum example, both of the obtained equations match the  exact system (\ref{pend}) with a good precision. On the other hand, the Lorenz system is a more challenging system. Hence, by the first-order approximation of $\dot x$ with $h=5\mathrm{ms}$, only an approximation of the dynamics can be obtained, while the exact system (\ref{lorenz}) is identified if we use the exact $\dot x$ . As shown in Fig. ~\ref{fig:lorenz_state} and Fig. ~\ref{fig:lorenz_vpe}, although  the model obtained for Lorenz system by using the approximate state variables does not closely match the exact dynamics, the obtained controller can successfully solve the regulation problem as long as the prediction errors remain bounded.    
\begin{table*}[]
	\normalsize
	
	\begin{center}
		\caption{The system dynamics and the corresponding value function obtained by the proposed method, where the exact and the approximated derivatives of the state variables are used in different scenarios  }\label{tb:pend_lorenz}
		\renewcommand{\arraystretch}{1.2}%
		\begin{tabular}{c||c}
			
			\hline
			\cellcolor{blue!5} Exact $\dot{x}$ & \cellcolor{blue!5} $\dot{x}\approx (x_{k+1}-x_k)/h, h=5ms$ \\
			\hline
			\hline
			\multicolumn{2}{|c|}{\cellcolor{green!5}Pendulum ($\Phi=\{1,x,\sin x\}$)}  \\\hline \hline
			$\begin{aligned} 
			&\dot{x_1}=-1.000x_2 \ \nonumber \\ 
			&\dot{x_2}=-1.000x_2  -19.600\sin(x_1)  + 40.000u
			\end{aligned}$ & $\begin{aligned} 
			&\dot{x_1}=-1.011x_2 \ \nonumber \\ 
			&\dot{x_2}=-0.995x_2  -19.665sin(x_1)  + 40.098u
			\end{aligned}$  \\ \hline
			$\begin{aligned} 
			V(x)=&1.974x_1^2 -0.058x_2x_1  +0.036x_2^2\nonumber\\
			&-2.2\sin(x_1)x_1 -0.077\sin(x_1)x_2 \nonumber\\
			&+1.548\sin^2(x_1)
			\end{aligned}$ & $\begin{aligned} 
			V(x) =&   2.049x_1^2 -0.058x_2x_1  +0.036x_2^2\nonumber\\ &-2.371\sin(x_1)x_1 -0.077\sin(x_1)x_2 \nonumber\\ &+1.630\sin^2(x_1)
			\end{aligned}$ \\ \hline \hline
			\multicolumn{2}{|c|}{\cellcolor{green!5}Chaotic Lorenz System ($\Phi=\{1,x,x^2,x^3,x_ix_j\},\quad i,j \in \{1,\dots,n\},i\neq j$)} \\\hline\hline
			$\begin{aligned} \\
			&\dot{x_1}=-10.000x_1   +10.000x_2    +1.000u  \nonumber\\
			&\dot{x_2}= 28.000x_1   -1.000x_2   -1.000x_1x_3 \nonumber\\
			&\dot{x_3}=-2.667x_3    +1.000x_1x_2 
			\end{aligned}$ & $\begin{aligned}
			&\dot{x_1}=-10.070x_1   +9.973x_2      +0.989u\nonumber\\
			&\dot{x_2}= 0.993x_1    -0.997x_2    +8.483x_3   -1.000x_1x_3\nonumber\\
			&\dot{x_3}=-8.483x_1   -8.483x_2   -2.666x_3    +1.000x_1x_2
			\end{aligned}$ \\
			\hline
			$\begin{aligned}
			V(x) =&  30.377x_1^2 +48.939x_2x_1 +25.311x_2^2\nonumber\\
			&+1.500x_3^2  -1.873x_1x_2x_3  +4.719x_1^2x_2^2\nonumber\\ &-3.291x_1^2x_3+1.469x_1x_3^2-0.012x_1^2x_2x_3
			\end{aligned}$ & $\begin{aligned}
			
			V(x) =& 11.193x_1^2 + 8.389x_2x_1 +42.855x_2^2\nonumber\\
			&-20.950x_3x_1 +28.441x_3x_2 +32.045x_3^2\nonumber\\ &-1.899x_1^2x_2 -4.777x_1x_2^2  -0.456x_1x_2x_3\nonumber\\  &+5.064x_1^2x_2^2  +2.953x_1^2x_3 -8.168x_1x_3^2 \nonumber\\
			& -2.633x_1^2x_3x_2  +1.353x_1^2x_3^2
			\end{aligned}$  \\ \hline
		\end{tabular}
	\end{center}
\end{table*}

\subsubsection{\textbf{Example 1}. (Pendulum)}
The state space description of the system is given as 
\begin{align}\label{pend}
&\dot{x_1}=-x_2,\ \nonumber \\ 
&\dot{x_2}=-\frac{g}{l}\sin(x_1)-\frac{k}{m}x_2+\frac{1}{ml^2}u,
\end{align}
where $m=0.1kg$, $l=0.5m$, $k=0.1$, and $g=9.8m/s^2$. The performance criteria are defined by the choices of $Q=\text{diag}([1,1])$, $R=2$.

\textbf{Objective}: The system is regulated to the unstable equilibrium point given by $x_\text{ref}\equiv0$.

 In Fig.~\ref{fig:pend}(b) and Fig.~\ref{fig:pend1}.(b) state and control responses for two runs of the learned system are illustrated starting from two different initial conditions. In table~\ref{tb:pend_lorenz}, the learned dynamics and value function are listed for the exact and the approximated $\dot x$.

 \subsubsection{\textbf{Example 2}. (Chaotic Lorenz System)} The system dynamics are defined by
\begin{align}\label{lorenz}
&\dot{x_1}=\sigma(x_2-x_1)+u,\nonumber\\
&\dot{x_2}=-x_2+x_1(\rho-x_3),\nonumber\\
&\dot{x_3}=x_1x_2-\beta x_3,
\end{align}

where $\sigma=10$, $\rho=28$, and $\beta=8/3$. Furthermore, we set the performance criteria to $Q=\text{diag}([160,160,12])$, $R=1$. This system has two unstable equilibrium points $(\pm\sqrt{72},\pm\sqrt{72},27)$, where the trajectories of the system oscillate around these points. 

\textbf{Objective}: By randomly setting $x_0\in [-40,40]^3$, we regulate the system to the unstable equilibrium $(-\sqrt{72},\\-\sqrt{72},27)$.

\begin{figure}[]
	\begin{center}
		\hspace{5mm}
		\includegraphics[width=7.3cm]{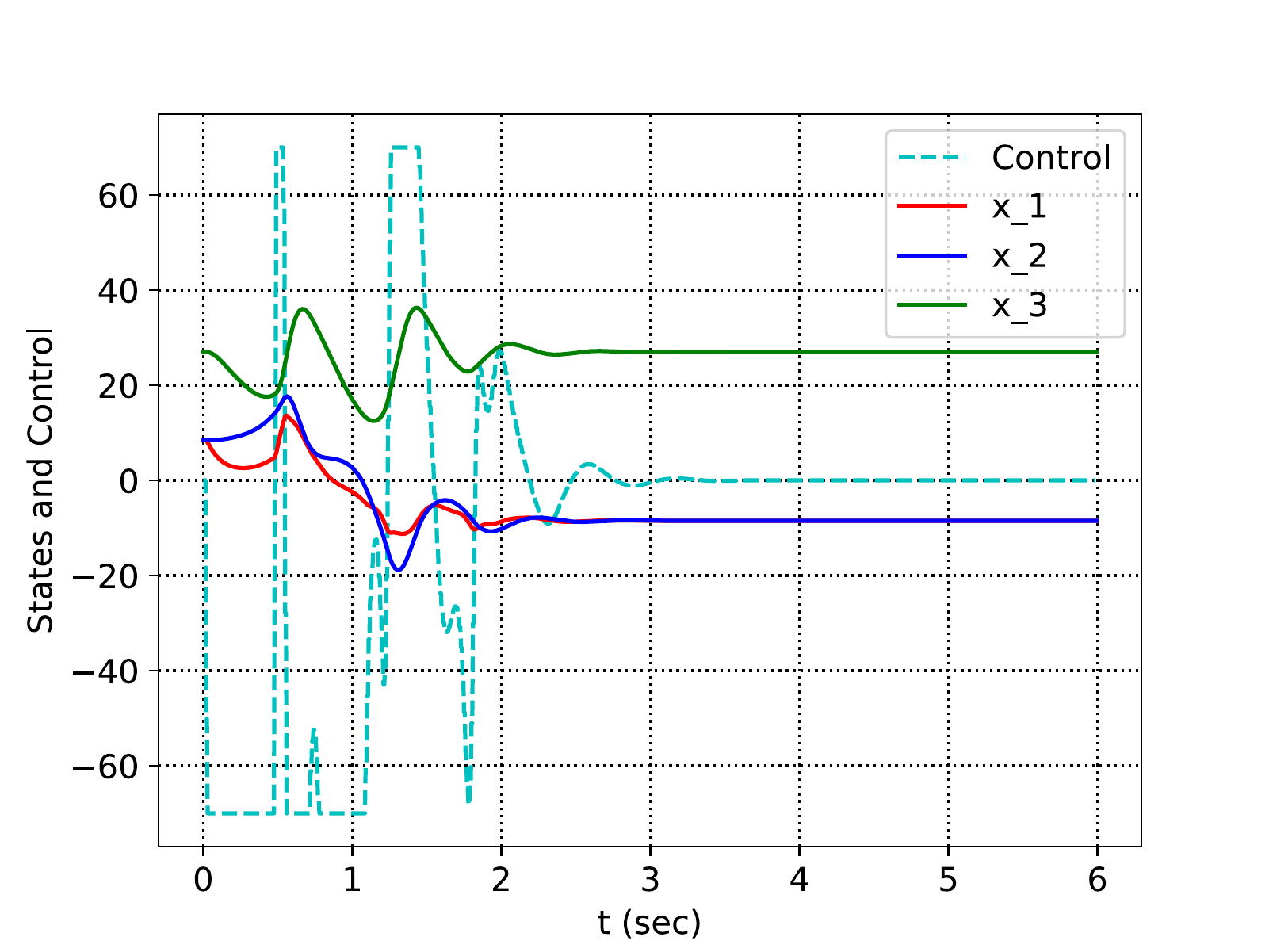}    %
		\caption{Responses of the Lorenz system while learning by using the approximated state derivatives as in Table \ref{tb:pend_lorenz}, where starting from one equilibrium point, we regulated the system to another unstable equilibrium.} 
		\label{fig:lorenz_state}
	\end{center}
\end{figure}

\begin{figure}[]
	\begin{center}
		\includegraphics[width=7cm]{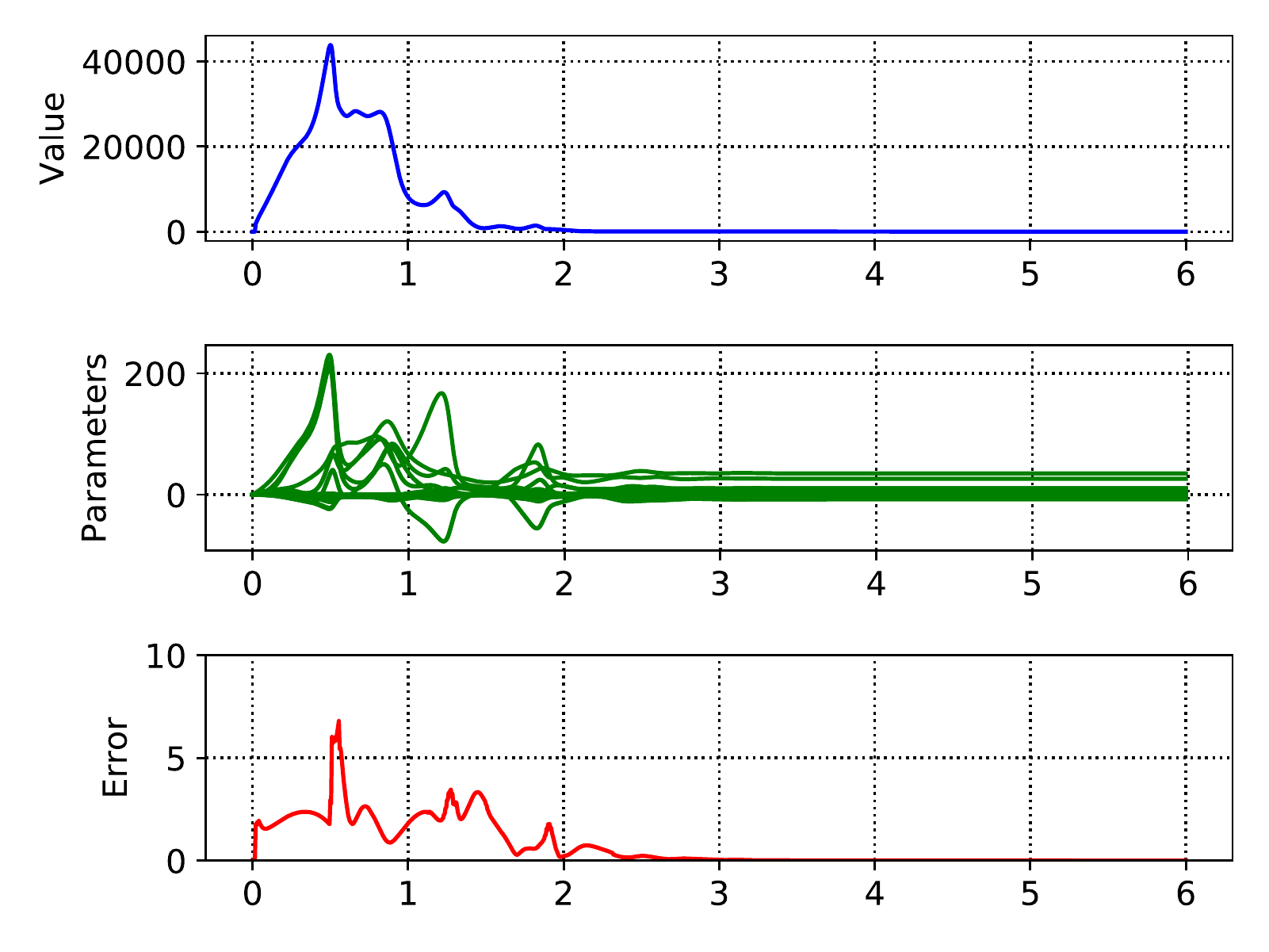}    %
		\caption{The value, components of $P$, and prediction error corresponding to Fig.~\ref{fig:lorenz_state}, respectively.} 
		\label{fig:lorenz_vpe}
	\end{center}
\end{figure}

\subsection{Systems to Be Approximated by a Given Set of Basis}
In what follows, we apply the presented learning scheme on two benchmark examples. Unlike the previous examples, the dynamics of these systems includes some rational terms that cannot be written in terms of some basis functions, however, an approximation can be obtained locally that is shown to be sufficient to successfully solve the regulation problem, as shown in Fig.~\ref{fig:cartpole_state}-\ref{fig:double_vpe}.  
\begin{figure}[h]
	\begin{center}
		\includegraphics[width=6cm]{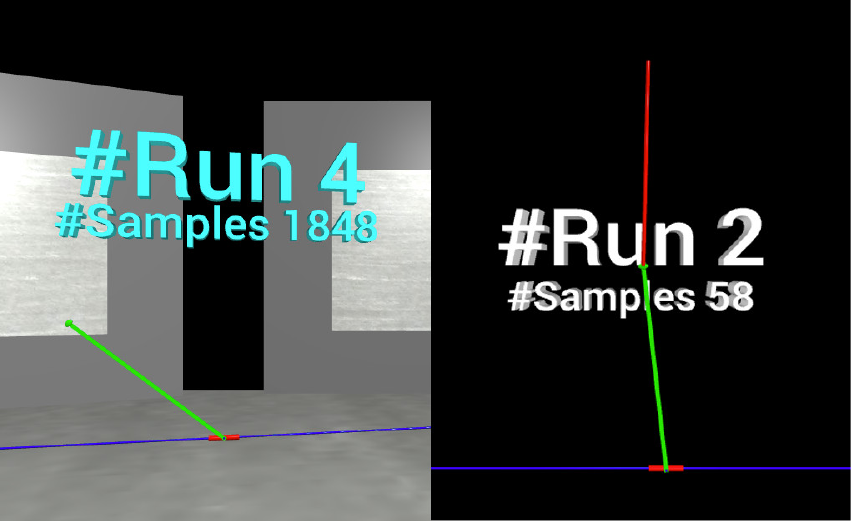}    %
		\caption{A view of the graphical simulations of the benchmark cartpole and double inverted pendulum examples. The video can be accessed in:\hspace{20mm} \protect\url{https://youtu.be/-j0vaHE9MZY} .} 
		\label{fig:3D}
	\end{center}
\end{figure}

Moreover, as shown in Fig.~\ref{fig:3D}, a video of the graphical simulation of the following benchmark examples is included. 
\subsubsection{\textbf{Example 3}. (Cartpole Swing up)}
The dynamics are given as
\begin{align}\label{cartpole}
&\dot{x_1}=x_2, \nonumber\\
&\dot{x_2}=\nonumber\\
&\quad\frac{-u\cos(x_1)-mLx_2^2\sin(x_1)\cos(x_1)+(M+m)g\sin(x_1)}{L(M+m\sin(x_1)^2)},\nonumber\\
&\dot{x_3}=x_4,\nonumber\\
&\dot{x_4}=\frac{u+m\sin(x_1)(Lx_2^2-g\cos(x_1))}{M+m\sin(x_1)^2},
\end{align}
where the state vector is composed of the angle of the pendulum from upright position, the angular velocity, and the position and velocity of the cart, with  
$m=0.1kg$, $M=1kg$, $L=0.8m$, and $g=9.8m/s^2$. Moreover, we choose $Q=\text{diag}([60,1.5,180,45])$, $R=1$. 
\vspace{-1mm}

\textbf{Objective}: By starting from some initial angles close to the stable angle of the pendulum ($\pm \pi$), the cart swings up the pendulum to reach to and stay at the unstable state given as $x_\text{ref}\equiv 0$.

By running the learning scheme, an approximation of the system is identified as  
\begin{align}
&\dot{x}_1=1.000x_2, \nonumber\\
&\dot{x}_2=12.934\sin(x_1)    +0.230\sin(x_3) -1.234\cos(x_1)u,\nonumber\\
&\dot{x}_3= 0.995x_4,\nonumber\\
&\dot{x}_4= 0.926\sin(x_1)    +0.953u,
\end{align}
where $\Phi=\{1,x,x^2,x^3,\sin x,\cos x\}$. 
Moreover, considering the assumed bases, we obtained the optimal value function as below.
\vspace{-3mm}
\begin{figure}[h!]
	\begin{center}
		\hspace{5mm}
		\includegraphics[width=7.5cm]{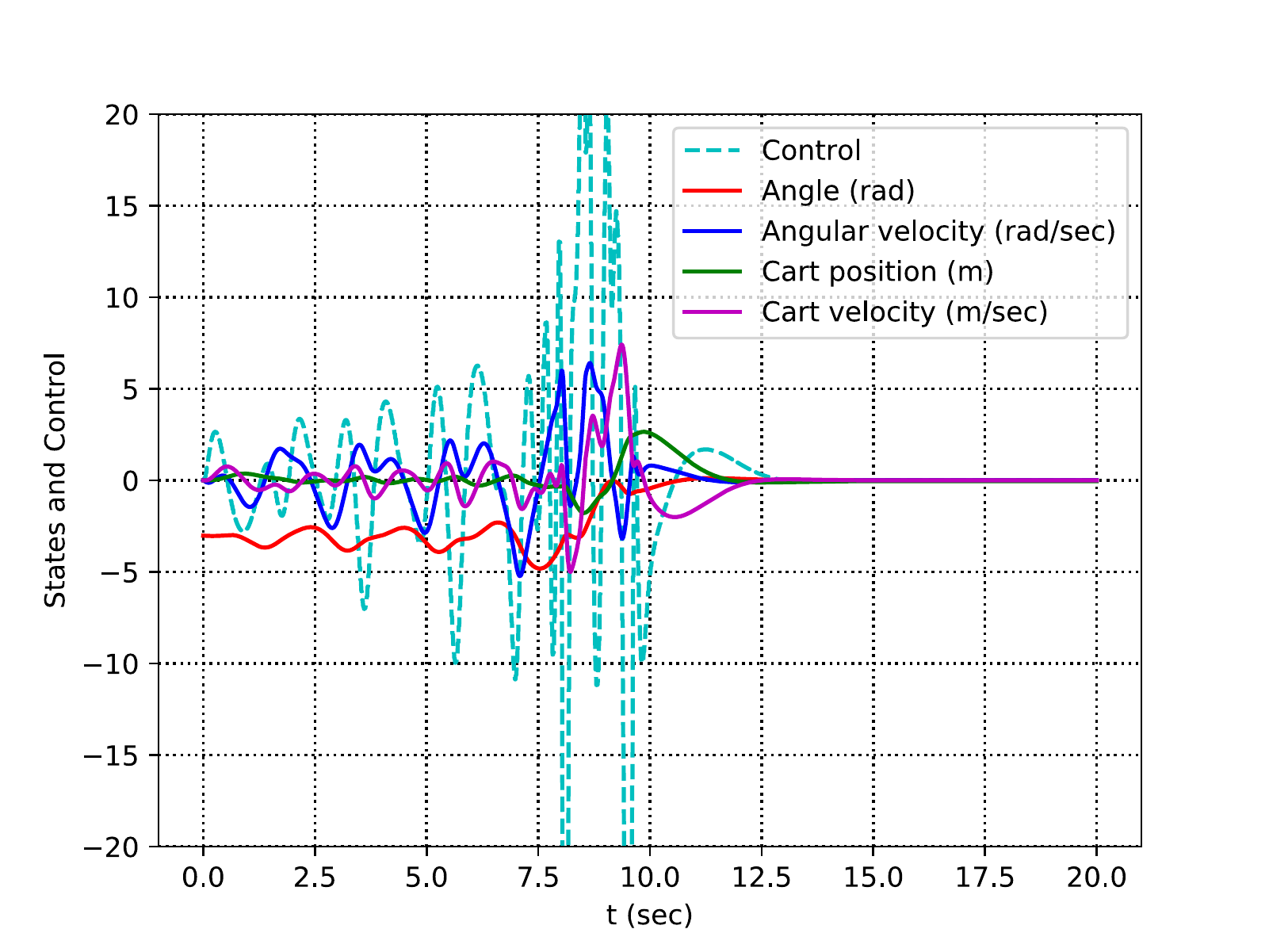}    %
		\caption{Responses of the cartpole system while learning by using the approximated state derivatives.}
		\label{fig:cartpole_state}
	\end{center}
\end{figure}
\begin{figure}[h!]
	\begin{center}
		\includegraphics[width=7.1cm]{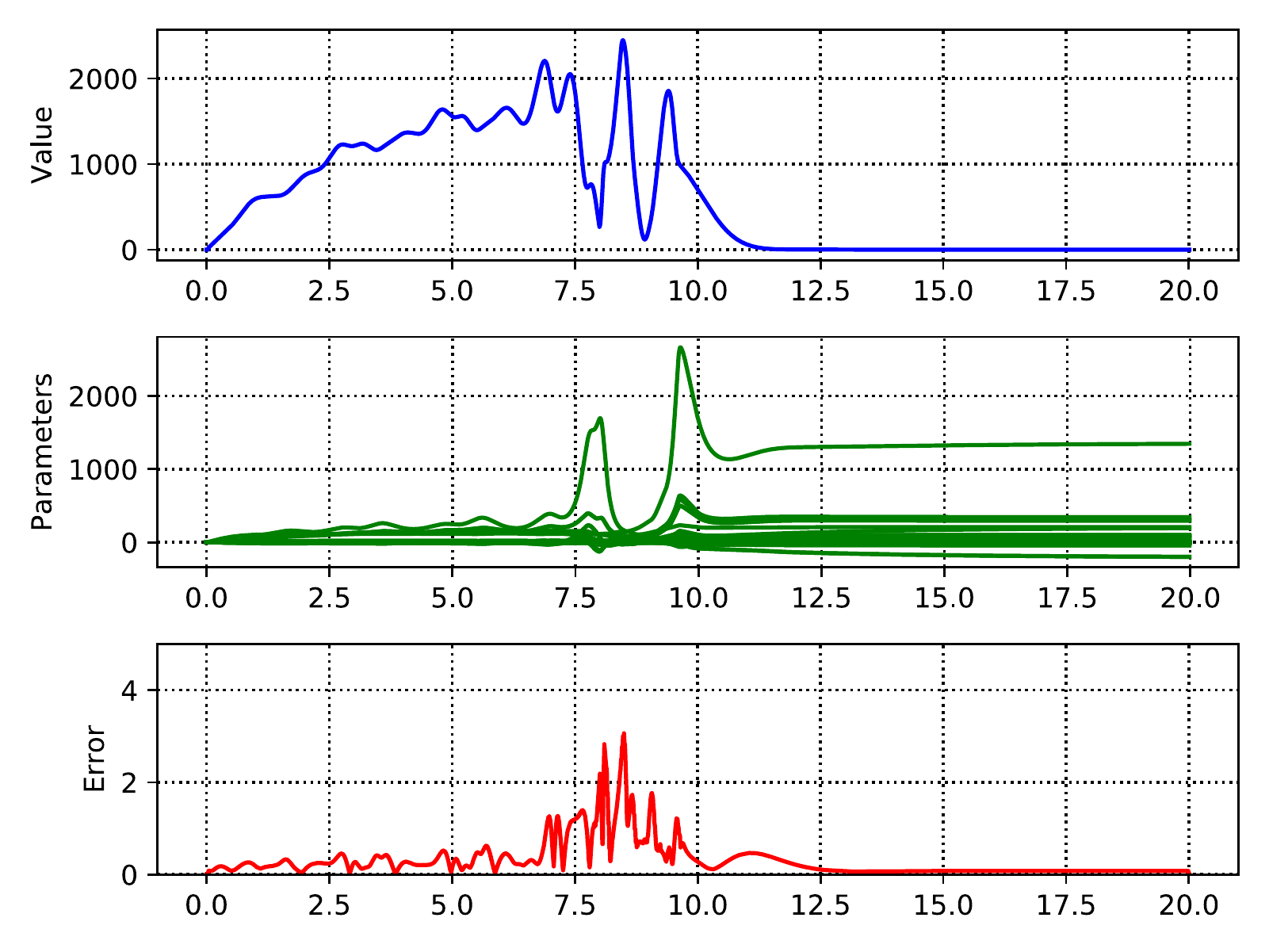}    %
		\caption{The value, components of $P$, and prediction error corresponding to Fig.~\ref{fig:cartpole_state}, respectively.} 
		\label{fig:cartpole_vpe}
	\end{center}
\end{figure}
\begin{align}
&V(x) =  59.712x_1^2  +9.855x_2x_1+134.855x_2^2+9.587x_3x_1 \nonumber\\ &+241.295x_3x_2+223.389x_3^2+4.418x_4x_1+222.022x_4x_2\nonumber\\ 
&+226.646x_4x_3+100.417x_4^2-63.050\sin(x_1)x_1\nonumber\\
&+1098.765\sin(x_1)x_2+2294.259\sin^2(x_1)\nonumber\\
&+984.786\sin(x_1)x_3+909.030\sin(x_1)x_4-1.712\sin(x_3)x_1\nonumber\\
 &+18.102\sin(x_3)x_2+15.812\sin(x_3)x_3+0.806\sin^2(x_3)\nonumber\\ &+75.231\sin(x_3)\sin(x_1)+15.072\sin(x_3)x_4.  
\end{align}

\subsubsection{\textbf{Example 4}. (Double Inverted Pendulum on a Cart)} By defining $y:=[q\quad \theta_1\quad \theta_2]^T$ to be a vector of the cart position and angles of the double pendulum from the top equilibrium point, the system dynamics can be written in the following form.
\begin{align} \label{double_pend}
\dot{x}=\begin{bmatrix}
&\dot{y}\\
&M^{-1}f(y,\dot{y})
\end{bmatrix},
\end{align}
where $M=$
\begin{align}
&\begin{bmatrix}
&m+m_1 +m_2, \hspace{2mm}l_1(m_1 +m_2)\cos(\theta_1), \hspace{2mm}m_2l_2 \cos(\theta_2)\\
&l_1(m_1 +m_2)\cos(\theta_1) , \hspace{0.5mm}l_1^2(m_1 +m_2), \hspace{0.5mm}l_1l_2m_2 \cos(\theta_1-\theta_2)\\
&l_2m_2 \cos(\theta_2), \hspace{3mm}l_1l_2m_2 \cos(\theta_1-\theta_2) , \hspace{8mm}l_2^2m_2
\end{bmatrix}, \nonumber
\end{align} 
\begin{align}
&f(y,\dot{y})= \nonumber\\
&\begin{bmatrix}
&l_1(m_1 +m_2)\dot\theta_1^2\sin(\theta_1)+m_2l_2 \dot\theta_2^2\sin(\theta_2)-d_1\dot{q}+u\\
&-l_1l_2m_2\dot\theta_2^2\sin(\theta_1-\theta_2)+g(m_1+m_2)l_1\sin(\theta_1)-d_2\theta_1\\
&l_1l_2m_2\dot\theta_1^2\sin(\theta_1-\theta_2)+gl_2m_2\sin(\theta_2)-d_3\theta_2
\end{bmatrix},\nonumber
\end{align}
$m=6kg$, $m_1=3kg$, $m_2=1kg$, $l_1=1m$, $l_2=2m$,$d_1=10$, $d_2=1$, and $d_3=0.5$. 
\vspace{-1mm}

\textbf{Objective}: We run the system from random angles around the top unstable equilibria of the pendulums given by $\theta_1=0$ and $\theta_2=0$, where the controller has to learn to regulate the system to $x_\text{ref}\equiv0$. 

We choose the bases as $\Phi=\{1,x,x^2\}$. Moreover the performance criteria is given by  $Q=\text{diag}([15,15,15,1,1,1])$, $R=1$. A sample of the obtained approximate dynamics is
\begin{align}\label{double_id}
&\dot{x}_1 = 0.998x_4,\quad \dot{x}_2 = 0.997x_5,\quad \dot{x}_3 = 0.996x_6, \nonumber\\
&\dot{x}_4 =    0.238x_1   -4.569x_2    1.245x_3   -1.891x_4   -0.908x_6   \nonumber\\
&\qquad -0.105x_2^2 +5.0131u   -2.824x_2^2u, \nonumber\\
&\dot{x}_5 =   16.718x_2   -2.328x_3    +1.558x_4 -0.598x_5  +0.130x_6 \nonumber\\
& \qquad-0.114x_2^2   -4.9911u    5.777x_2^2u   -0.690x_3^2u, \nonumber\\
&\dot{x}_6 =    0.123x_1   -6.721x_2    +9.032x_3    +0.191x_5   -0.358x_6    \nonumber\\
&\qquad +0.969x_3^2+0.184x_6^2   -1.898x_2^2u    +1.431x_3^2u.
\end{align}
It should be noted that, because of the random initial conditions and different samples in the database, a different approximation of the system may be obtained in any learning procedure. Furthermore, considering the dimension of the system and the number of terms in the identified system (\ref{double_id}), the obtained value function includes many terms of polynomials as expected. Therefore, for the sake of brevity, the obtained optimal value function is omitted.

\section{Conclusion}

Considering the online model-based regulation problem, the structured dynamics helped us in analytically computing an iterative update rule to improve the optimal value function according to the latest update on the identified system. Based on the computational complexity and the performance observed in the numerical and graphical simulations, we showed some potential opportunities in employing the SOL algorithm as an online model-based learning technique. Our future research will follow on the stability analysis and further applications of this approach.              
\vspace{-6mm}
\begin{figure}[h]
	\begin{center}
		\hspace{5mm}
		\includegraphics[width=7.2cm]{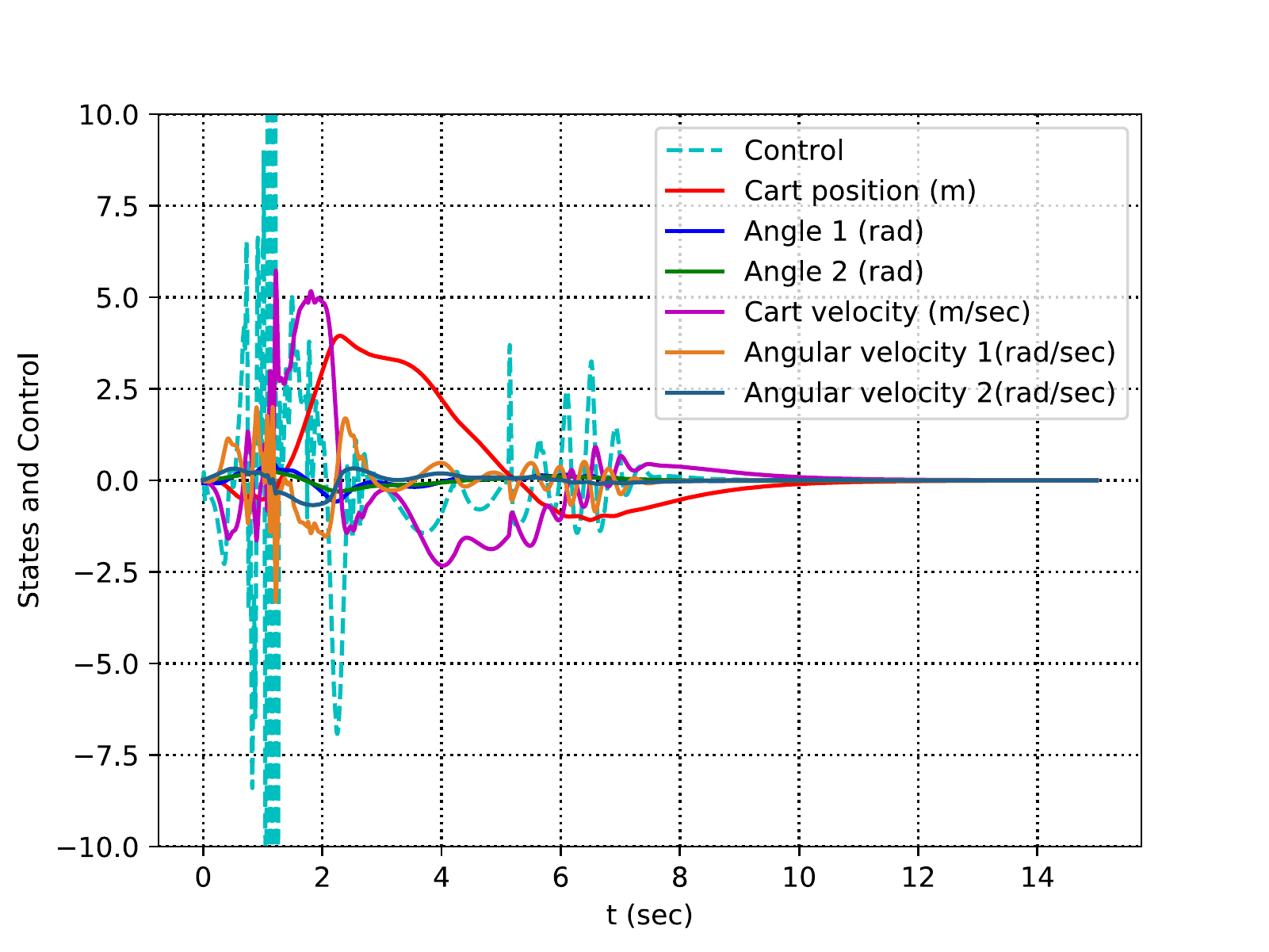} \vspace{-4mm}   %
		\caption{Responses of the double-inverted pendulum system while learning by using the approximated state derivatives.} 
		\label{fig:double_state}
	\end{center}
\end{figure}

\begin{figure}[h]
	\begin{center}
		\includegraphics[width=6.8cm]{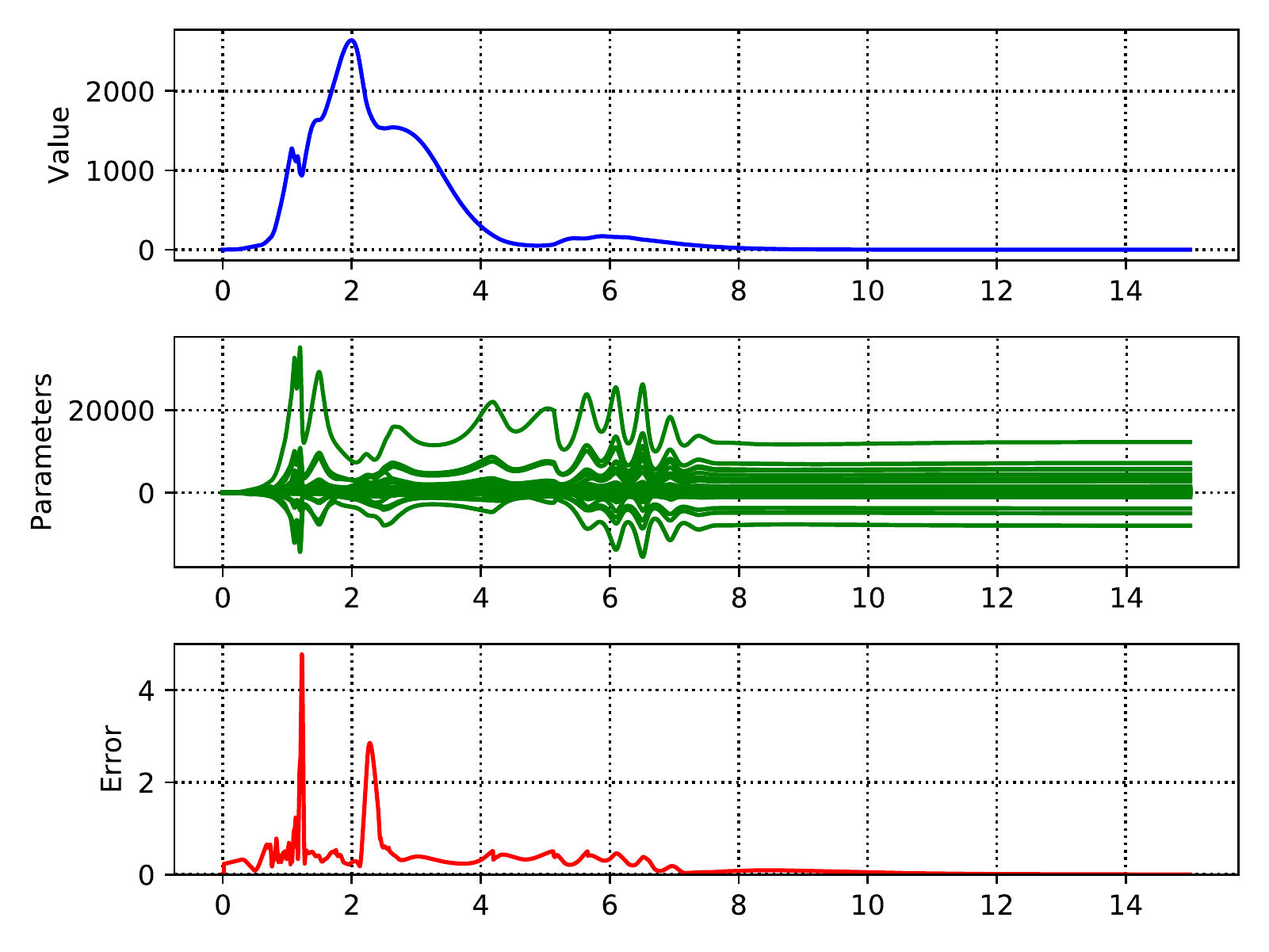}    %
		\caption{The value, components of $P$, and prediction error corresponding to Fig.~\ref{fig:double_state}, respectively.} 
		\label{fig:double_vpe}
	\end{center}
\end{figure}

\bibliography{ifacconf}             %

\begin{thebibliography}{22}
\providecommand{\natexlab}[1]{#1}
\providecommand{\url}[1]{\texttt{#1}}
\providecommand{\urlprefix}{URL }
\expandafter\ifx\csname urlstyle\endcsname\relax
  \providecommand{\doi}[1]{doi:\discretionary{}{}{}#1}\else
  \providecommand{\doi}{doi:\discretionary{}{}{}\begingroup
  \urlstyle{rm}\Url}\fi

\bibitem[{Atkeson and Santamaria(1997)}]{SOL:atkeson1997comparison}
Atkeson, C.G. and Santamaria, J.C. (1997).
\newblock A comparison of direct and model-based reinforcement learning.
\newblock In \emph{Proceedings of International Conference on Robotics and
  Automation}, volume~4, 3557--3564. IEEE.

\bibitem[{Barto et~al.(1983)Barto, Sutton, and
  Anderson}]{SOL:barto1983neuronlike}
Barto, A.G., Sutton, R.S., and Anderson, C.W. (1983).
\newblock Neuronlike adaptive elements that can solve difficult learning
  control problems.
\newblock \emph{IEEE transactions on systems, man, and cybernetics}, (5),
  834--846.

\bibitem[{Benosman(2018)}]{SOL:benosman2018model}
Benosman, M. (2018).
\newblock Model-based vs data-driven adaptive control: An overview.
\newblock \emph{International Journal of Adaptive Control and Signal
  Processing}, 32(5), 753--776.

\bibitem[{Bhasin et~al.(2013)Bhasin, Kamalapurkar, Johnson, Vamvoudakis, Lewis,
  and Dixon}]{SOL:bhasin2013novel}
Bhasin, S., Kamalapurkar, R., Johnson, M., Vamvoudakis, K.G., Lewis, F.L., and
  Dixon, W.E. (2013).
\newblock A novel actor--critic--identifier architecture for approximate
  optimal control of uncertain nonlinear systems.
\newblock \emph{Automatica}, 49(1), 82--92.

\bibitem[{Brunton et~al.(2016)Brunton, Proctor, and
  Kutz}]{SOL:brunton2016discovering}
Brunton, S.L., Proctor, J.L., and Kutz, J.N. (2016).
\newblock Discovering governing equations from data by sparse identification of
  nonlinear dynamical systems.
\newblock \emph{Proceedings of the National Academy of Sciences}, 113(15),
  3932--3937.

\bibitem[{Busoniu et~al.(2017)Busoniu, Babuska, De~Schutter, and
  Ernst}]{SOL:2017reinforcement}
Busoniu, L., Babuska, R., De~Schutter, B., and Ernst, D. (2017).
\newblock \emph{Reinforcement learning and dynamic programming using function
  approximators}.
\newblock CRC press.

\bibitem[{Duan et~al.(2016)Duan, Chen, Houthooft, Schulman, and
  Abbeel}]{SOL:duan2016benchmarking}
Duan, Y., Chen, X., Houthooft, R., Schulman, J., and Abbeel, P. (2016).
\newblock Benchmarking deep reinforcement learning for continuous control.
\newblock In \emph{International Conference on Machine Learning}, 1329--1338.

\bibitem[{Kaiser et~al.(2018)Kaiser, Kutz, and Brunton}]{SOL:kaiser2018sparse}
Kaiser, E., Kutz, J.N., and Brunton, S.L. (2018).
\newblock Sparse identification of nonlinear dynamics for model predictive
  control in the low-data limit.
\newblock \emph{Proceedings of the Royal Society A}, 474(2219), 20180335.

\bibitem[{Kamalapurkar et~al.(2016{\natexlab{a}})Kamalapurkar, Rosenfeld, and
  Dixon}]{SOL:kamalapurkar2016efficient}
Kamalapurkar, R., Rosenfeld, J.A., and Dixon, W.E. (2016{\natexlab{a}}).
\newblock Efficient model-based reinforcement learning for approximate online
  optimal control.
\newblock \emph{Automatica}, 74, 247--258.

\bibitem[{Kamalapurkar et~al.(2016{\natexlab{b}})Kamalapurkar, Walters, and
  Dixon}]{SOL:kamalapurkar2016model}
Kamalapurkar, R., Walters, P., and Dixon, W.E. (2016{\natexlab{b}}).
\newblock Model-based reinforcement learning for approximate optimal
  regulation.
\newblock \emph{Automatica (Journal of IFAC)}, 64(C), 94--104.

\bibitem[{Kamalapurkar et~al.(2018)Kamalapurkar, Walters, Rosenfeld, and
  Dixon}]{SOL:kamalapurkar2018model}
Kamalapurkar, R., Walters, P., Rosenfeld, J., and Dixon, W. (2018).
\newblock Model-based reinforcement learning for approximate optimal control.
\newblock In \emph{Reinforcement Learning for Optimal Feedback Control},
  99--148. Springer.

\bibitem[{Kivinen et~al.(2004)Kivinen, Smola, and
  Williamson}]{SOL:kivinen2004online}
Kivinen, J., Smola, A.J., and Williamson, R.C. (2004).
\newblock Online learning with kernels.
\newblock \emph{IEEE Transactions on Signal Processing}, 52(8), 2165--2176.

\bibitem[{Lewis and Vrabie(2009)}]{SOL:lewis2009reinforcement}
Lewis, F.L. and Vrabie, D. (2009).
\newblock Reinforcement learning and adaptive dynamic programming for feedback
  control.
\newblock \emph{IEEE Circuits and Systems Magazine}, 9(3), 32--50.

\bibitem[{Modares et~al.(2014)Modares, Lewis, and
  Naghibi-Sistani}]{SOL:modares2014integral}
Modares, H., Lewis, F.L., and Naghibi-Sistani, M.B. (2014).
\newblock Integral reinforcement learning and experience replay for adaptive
  optimal control of partially-unknown constrained-input continuous-time
  systems.
\newblock \emph{Automatica}, 50(1), 193--202.

\bibitem[{Polydoros and Nalpantidis(2017)}]{SOL:polydoros2017survey}
Polydoros, A.S. and Nalpantidis, L. (2017).
\newblock Survey of model-based reinforcement learning: Applications on
  robotics.
\newblock \emph{Journal of Intelligent \& Robotic Systems}, 86(2), 153--173.

\bibitem[{Powell(2004)}]{SOL:powell2004handbook}
Powell, W.B. (2004).
\newblock \emph{Handbook of learning and approximate dynamic programming},
  volume~2.
\newblock John Wiley \& Sons.

\bibitem[{Recht(2019)}]{SOL:recht2019tour}
Recht, B. (2019).
\newblock A tour of reinforcement learning: The view from continuous control.
\newblock \emph{Annual Review of Control, Robotics, and Autonomous Systems}, 2,
  253--279.

\bibitem[{Scherer et~al.(2000)Scherer, Dubois, and
  Sherwood}]{SOL:scherer2000vpython}
Scherer, D., Dubois, P., and Sherwood, B. (2000).
\newblock Vpython: 3d interactive scientific graphics for students.
\newblock \emph{Computing in Science \& Engineering}, 2(5), 56--62.

\bibitem[{Sutton(1990)}]{SOL:sutton1990integrated}
Sutton, R.S. (1990).
\newblock Integrated architectures for learning, planning, and reacting based
  on approximating dynamic programming.
\newblock In \emph{Machine Learning Proceedings 1990}, 216--224. Elsevier.

\bibitem[{Vamvoudakis et~al.(2012)Vamvoudakis, Lewis, and
  Hudas}]{SOL:vamvoudakis2012multi}
Vamvoudakis, K.G., Lewis, F.L., and Hudas, G.R. (2012).
\newblock Multi-agent differential graphical games: Online adaptive learning
  solution for synchronization with optimality.
\newblock \emph{Automatica}, 48(8), 1598--1611.

\bibitem[{Van~Vaerenbergh and Santamar{\'\i}a(2014)}]{SOL:van2014online}
Van~Vaerenbergh, S. and Santamar{\'\i}a, I. (2014).
\newblock Online regression with kernels.
\newblock In \emph{Regularization, Optimization, Kernels, and Support Vector
  Machines}, 495--521. Chapman and Hall/CRC.

\bibitem[{Zhang et~al.(2011)Zhang, Cui, Zhang, and Luo}]{SOL:zhang2011data}
Zhang, H., Cui, L., Zhang, X., and Luo, Y. (2011).
\newblock Data-driven robust approximate optimal tracking control for unknown
  general nonlinear systems using adaptive dynamic programming method.
\newblock \emph{IEEE Transactions on Neural Networks}, 22(12), 2226--2236.

\end{thebibliography}

\end{document}